\documentclass{amsart}
\usepackage[T1]{fontenc}
\usepackage{latexsym,amssymb,amsmath}
\usepackage{amsthm}
\usepackage{longtable}
\usepackage{epsfig}
\usepackage{hhline}
\usepackage{epic}

    \newcommand{\Res}{\operatorname{Res}}
    \newcommand{\Ind}{\operatorname{Ind}}
    
    \newcommand{\Cen}{\operatorname{C}}
    \newcommand{\MP}{\mathcal{MP}}

    \newcommand{\cyc}[1]{\langle\,#1\,\rangle}
    
    \newcommand{\Irr}{\operatorname{Irr}}
    
    \newcommand{\sym}{\mathfrak{S}}

    \newcommand{\mm}{\boldsymbol{\mu}}

    \newcommand{\Z}{\mathbb{Z}}
    \newcommand{\A}{\mathfrak A}
    \newcommand{\B}{\cal B}
    \newcommand{\bt}{\boxtimes}

\newcommand{\cal}[1]{\mathcal{#1}}

\newtheorem{theorem}{Theorem}[section] % 1st argument is your name for it
\newtheorem{lemma}[theorem]{Lemma}     % 2nd argument is what is printed

\newtheorem{proposition}[theorem]{Proposition}

\newtheorem{remark}[theorem]{Remark}

\title[A $2$-basic set of the alternating group]% end with percent
    {A $2$-basic set of the alternating group}

\author{Olivier Brunat}

\address{Ruhr-Universit\"at Bochum\\
Fakult\"at f\"ur Mathematik\\
Raum NA 2/33\\
D-44780 Bochum\\}
\email{Olivier.Brunat@ruhr-uni-bochum.de}

\author{Jean-Baptiste Gramain}

\address{\'Ecole Polytechnique F\'ed\'erale de Lausanne\\
Institut de g\'eom\'etrie, alg\`ebre et topologie\\
B\^atiment de Chimie (BCH)\\
CH-1015 Lausanne\\}
\email{jean-baptiste.gramain@epfl.ch}

\begin{document}
\maketitle
\begin{abstract}In this note, we construct a $2$-basic set of the
alternating group $\A_n$. To do this, we construct a $2$-basic set of
the
symmetric group $\sym_n$ with an additional property, such that its
restriction to $\A_n$ is a $2$-basic set.
We adapt here a method developed in~\cite{BrGr} for the case when
the characteristic is odd. One of the main tools is the generalized
perfect isometries defined by K\"ulshammer, Olsson and Robinson
in~\cite{KOR}.
\end{abstract}
\section{Introduction}
This note is concerned with the existence problem of basic sets
for finite groups. Let $G$ be a finite group and $p$ be a prime.
A $p$-basic set of $G$ is a subset $B$ of the set $\Irr(G)$ of
irreducible complex characters of $G$, such that the restrictions
$B^{p\textrm{-reg}}=\{\chi^{p\textrm{-reg}}, \, \chi \in B \}$
to $p$-regular elements
(i.e., the elements of $G$ with order prime to $p$) of the characters
in $B$ form a $\Z$-basis of the ring of $p$-Brauer characters of $G$.
Gerhard Hiss has conjectured that every finite group has a $p$-basic
set. It is actually proved for some finite groups (see for example
\cite{GH}, \cite{GeckBS}, \cite{GeckBS3} and \cite{Br1}), but it remains
open and a difficult question in general.

Recently, the authors proved in~\cite{BrGr} that the alternating group
$\A_n$ has
a $p$-basic set for any odd prime $p$. In this note, we complete this work
by constructing a $2$-basic set of $\A_n$.

Recall that the irreducible characters of the symmetric
group $\sym_n$ are naturally labelled by the partitions of
$n$~\cite[2.1.11]{James-Kerber}. We denote by ${\cal P}$ (respectively
${\cal P_n}$) the set of all partitions of all integers (respectively of
$n$). For any $\lambda \in {\cal P}$, we write $| \lambda |$ for the size
of $\lambda$. For any partition $\lambda$ of
$n$ (written $\lambda\vdash n$ in the following), we denote by
$\chi_{\lambda}$ the corresponding irreducible character of $\sym_n$ and
by
$\lambda^*$ the conjugate partition of $\lambda$. Then $\chi_{\lambda^*}=
\varepsilon \chi_{\lambda}$, where $\varepsilon$ is the signature of
$\sym_n$. A partition $\lambda$ is said to be \emph{self-conjugate} if
$\lambda = \lambda^*$.

For any $B \subseteq \Irr(\sym_n)$, we call \emph{restriction of $B$ to
$\A_n$} the subset $B_{\A_n}$ of $\Irr(\A_n)$ consisting of all the
irreducible constituents of the restrictions to $\A_n$ of the characters
in $B$. To prove the
existence of a $p$-basic set of $\A_n$ for odd prime $p$,
the approach in~\cite{BrGr} consists in constructing a $p$-basic
set $\B$ of $\sym_n$ satisfying two additional
properties which ensure that $\B_{\A_n}$ is a $p$-basic set of $\A_n$ (cf.
\cite[1.1]{BrGr}).

We will see that, when $p=2$, we can use a similar strategy, by requiring
the single additional property that $\B$ contains all characters labelled
by self-conjugate partitions; namely, we will prove

%\begin{enumerate}
%\item If $\chi_{\lambda}\in\B_{\emptyset}$, then
%$\chi_{\lambda^*}\in\B_{\emptyset}$.
%\item If $\lambda=\lambda^*$, then $\chi_{\lambda}\in\B_{\emptyset}$ if
%and only if $\overline{\lambda}$ has no parts divisible by $p$, where, if
%$\lambda=(\lambda_1,\ldots,\lambda_r) \vdash n$, then $\overline{\lambda}
%= (2\lambda_1-1,2\lambda_2-3,\ldots,1) \vdash n$.
%\end{enumerate}

%Suppose now that $p=2$ and that $B$ is a $2$-basic set of $\sym_n$.
%Then Property (1) is never satisfied; see Lemma~\ref{ressign}.
% Indeed, for every $\lambda\vdash
%n$ the restrictions of $\chi_{\lambda}$ and $\chi_{\lambda^*}$ on
%$2$-regular elements coincide. Moreover, if $\lambda$ is self-conjugate,
%then
%every part of $\overline{\lambda}$ is odd. Therefore, for $p=2$,
%Properties (1) and (2) reduce to the following condition:
%the set $B$
%contains all self-conjugate partitions. We will prove that this property
%is an additional condition assuring that $2$-basic sets of $\sym_n$
%restrict to $2$-basic sets of $\A_n$. We recall that the restriction
%$B_{\A_n}$ of a set $B\subseteq\Irr(\sym_n)$ to $\A_n$ is the set of the
%irreducible constituents of $\Res_{\A_n}^{\sym_n}(\chi)$ for $\chi\in
%B$. More precisely, we will prove

\begin{proposition}\label{resAn}
If $B$ is a $2$-basic set of $\sym_n$ containing every
$\chi_{\lambda}\in\Irr(\sym_n)$ with $\lambda=\lambda^*$, then the
restriction $B_{\A_n}$ of $B$ to $\A_n$ is a $2$-basic set of $\A_n$.
\end{proposition}

The aim of this paper is to adapt the methods of~\cite{BrGr}
(some of which heavily relie on $p$ being odd) to the case $p=2$,
in order to construct a $2$-basic set of $\sym_n$ satisfying
Proposition~\ref{resAn}.\\

%Let $\lambda=(\lambda_1,\ldots,\lambda_r)$ be a partition of $n$.
%We consider an abacus with $2$ runners and fix an
%origin. Using it, we can construct the $2$-quotient
%$\alpha_{\lambda}=(\lambda^{(1)},\lambda^{(2)})$ of a partition
%$\lambda\vdash n$.
%Note that we can choose the origin of the abacus such that,
%for any $n$ and any $\lambda \vdash n$, if
%$\alpha_{\lambda}=(\lambda^{(1)},\lambda^{(2)})$, then
%$\alpha_{\lambda^*}=(\lambda^{(2)^*},\lambda^{(1)^*})$ (see the proof
%of~\cite[Lemma 3.1]{BrGr}). From now on, we
%always assume this to be the case. Then our main result is

%We denote by $\alpha_{\lambda}=(\lambda^{(1)},\lambda^{(2)})$
%the $2$-quotient of $\lambda$. Note that this is only defined up
%to a choice of convention. However, it is proved in~\cite{BrGr}
%(cf Lemma 3.1) that, for a certain choice,
%we have that, for any $n$ and any $\lambda \vdash n$, if
%$\alpha_{\lambda}=(\lambda^{(1)},\lambda^{(2)})$, then
%$\alpha_{\lambda^*}=(\lambda^{(2)^*},\lambda^{(1)^*})$. From now on, we
%always assume this to be the case. Then our main result is

Let $\lambda=(\lambda_1,\ldots,\lambda_r)$ be a partition of $n$.
We denote by $\alpha_{\lambda}=(\lambda^{(1)},\lambda^{(2)})$
the $2$-quotient of $\lambda$. Note that this is only defined up
to a choice of convention. However, it is proved in~\cite{BrGr}
(cf. the proof of Lemma 3.1) that, for a certain choice,
we have that, for any $n$ and any $\lambda \vdash n$, if
$\alpha_{\lambda}=(\lambda^{(1)},\lambda^{(2)})$, then
$\alpha_{\lambda^*}=(\lambda^{(2)^*},\lambda^{(1)^*})$. From now on, we
always assume this to be the case. Then our main result is

\begin{theorem}\label{main}For every integer $m$, we set
$$\cal P_m'=\{\lambda=(\lambda_1,\ldots,\lambda_r)\vdash
m\,|\,\lambda_i\ \textrm{is even for }1\leq i\leq r\}.$$
Define
$$\Lambda=\{(\mu,\emptyset)\ |\ \mu\notin\cal
P_{|\mu|}'\}\cup\{(\mu,\mu^*), \, \mu \in {\cal P}\}.$$
Then the set $\cal B_{\Lambda}$ of irreducible characters
$\chi_{\lambda}$ of $\sym_n$ satisfying $\alpha_\lambda\in\Lambda$
is a $2$-basic set of $\sym_n$. Moreover the restriction $\cal
B_{\Lambda,\A_n}$ of $\cal B_{\Lambda}$ to $\A_n$ is a $2$-basic set of
$\A_n$.
\end{theorem}

Throughout this article, we will use the following notations and
conventions. For $G$ a finite group, $\Irr(G)$ denotes the set of
complex irreducible characters of $G$. For $\phi_1,\,\phi_2\in\Irr(G)$,
we denote by $\phi_1\otimes\phi_2$ the character defined by
$\phi_1\otimes\phi_2(g)=\phi_1(g)\phi_2(g)$. Let $H$ and $K$ be two
finite groups. For $\phi_H\in\Irr(H)$ and $\phi_K\in\Irr(K)$, we
define $\phi_H\bt\phi_K\in\Irr(H\times K)$ by
$$\phi_H\bt\phi_K(h,k)=\phi_H(h)\phi_K(k)\quad\textrm{for }h\in H,\,k\in
K.$$
For $c$ a conjugacy class of $G$ and $\chi$ a character of $G$, we
sometimes will use the notation $\chi(c)$ for $\chi(g)$ for $g\in c$.

The paper is organized as follows. In Section~\ref{part1}, we prove
Proposition~\ref{resAn}. In Section~\ref{part2}, we construct, for
each integer $w$, a $\Z$-basis of the ring $\Z\Irr(\sym_{2w})$ of
virtual characters of $\sym_{2w}$ that we will need in order to
prove Theorem~\ref{main}. Finally, in Section~\ref{part3}, we prove
Theorem~\ref{main}.

\section{Restriction to $\A_n$ of $2$-basic sets of $\sym_n$}\label{part1}

We first prove the following lemma.
\begin{lemma}
With the above notation, if $\lambda\vdash n$, then
$\chi_{\lambda}^{2\textrm{-reg}}=\chi_{\lambda^*}^{2\textrm{-reg}}$.
\label{ressign}
\end{lemma}
\begin{proof}
Let $\sigma$ be an element of $\sym_n$ with odd order. Write
$\sigma=c_1\cdots c_r$ as product of disjoint cycles. Since the
order of $\sigma$ is odd, each of the $c_i$'s must have odd length, so
that $\varepsilon (\sigma)=1$. Since $\chi_{\lambda^*}= \varepsilon
\chi_{\lambda}$, this yields the claim.
\end{proof}

\subsection{Irreducible characters of $\A_n$}\label{charAn}
The construction and the values of the irreducible characters of $\A_n$
are described in~\cite[2.5.13]{James-Kerber}. For the convenience of
the reader, we briefly recall how to parametrize them from those of
$\sym_n$. Take any $\lambda \vdash n$, and let
$$\rho_{\lambda}:=\Res_{\A_n}^{\sym_n}(\chi_{\lambda}).$$
If $\lambda\neq\lambda^*$, then $\rho_{\lambda}=\rho_{\lambda^*}$ is
irreducible, and $\Ind_{\A_n}^{\sym_n}(\rho_{\lambda})=\chi_{\lambda}+
\chi_{\lambda^*}$. Otherwise, $\rho_{\lambda}$ is the sum of two
irreducible characters of
$\A_n$, written $\rho_{\lambda,\pm}$, and chosen as follows. If
$\lambda=\lambda^*=(\lambda_1,\ldots,\lambda_r) \vdash n$, then we let
$\overline{\lambda}= (2\lambda_1-1,2\lambda_2-3,\ldots,1) \vdash n$.The
conjugacy
class of $\sym_n$ of cycle type $\overline{\lambda}$ consists of elements
of $\A_n$, and splits into two
classes $\overline{\lambda}_{\pm}$ of $\A_n$. Now, if $x\in\A_n$ doesn't
have cycle type $\overline{\lambda}$, then
$\rho_{\lambda,+}(x)=\rho_{\lambda,-}(x)$.
If $x_{\pm}\in\overline{\lambda}_{\pm}$, then
$\rho_{\lambda,\pm}(x_+)=s_{\lambda}\pm t_{\lambda}$ and
$\rho_{\lambda,\pm}(x_-)=s_{\lambda},\mp t_{\lambda}$, with $s_{\lambda}$
and $t_{\lambda}$ as described in~\cite[2.5.13]{James-Kerber}.
Furthermore,
$\Ind_{\A_n}^{\sym_n}(\rho_{\lambda,+})=\Ind_{\A_n}^{\sym_n}(\rho_{\lambda,-})=\chi_{\lambda}$.

\subsection{Proof of Proposition~\ref{resAn}}
The proof of Proposition~\ref{resAn} is based on the proof
of~\cite[5.2]{BrGr}.
%Note that in~\cite[5.2]{BrGr}, the basic
%set of $\sym_n$ is supposed to be $\varepsilon$-stable, where
%$\varepsilon$ is the sign character of $\sym_n$. As we will show,
%we can remove this assumption in our situation.
Suppose that $\B$ is a $2$-basic set of $\sym_n$ as in
Proposition~\ref{resAn} and consider its restriction $\B_{\A_n}$ to
$\A_n$.
%Denote by $\cal A$ the ring generated over $\Z$ by the restrictions
%$\chi^{2\textrm{-reg}}$ for $\chi\in \Irr(\A_n)$.
To prove that $\B_{\A_n}$ is a $2$-basic set of $\A_n$ we have
to show that $\B_{\A_n}^{2\textrm{-reg}}=\{\chi^{2\textrm{-reg}}|\chi\in
\B_{\A_n}\}$ is free and generates over $\Z$ the ring
$\Z\Irr(\A_n)^{2\textrm{-reg}}$.
%Note that in the proof of~\cite[5.2]{BrGr}, to prove that $\B_{\A_n}$
%is free, we do not need the $\varepsilon$-stable condition
%on $\B$.

We denote by $\cal S$ the set of self-conjugate partition of $n$ and by
$\cal T$ the set of partitions of $n$ such that $\lambda\in\cal T$ if
and only if $\chi_{\lambda}\in\B$. Put $\cal S'=\cal T\backslash \cal
S$.
With this notation, we have
$$\B_{\A_n}=\{\rho_{\lambda}\,|\,\lambda\in \cal
S'\}\cup\{\rho_{\lambda,\pm}\,|\,\lambda\in\cal S\}.$$
Suppose that there are integers $a_{\lambda}$ ($\lambda\in\cal S'$) and
$b_{\lambda,\pm}$ ($\lambda\in\cal S$) such that
\begin{equation}\label{eqlibre}
\sum_{\lambda\in\cal S'}a_{\lambda}\rho_{\lambda}^{2\textrm{-reg}}
+\sum_{\lambda\in\cal S}
(b_{\lambda,-}\rho_{\lambda,-}^{2\textrm{-reg}}+
b_{\lambda,+}\rho_{\lambda,+}^{2\textrm{-reg}})=0.
\end{equation}
For $\lambda\in\cal S'$, we have $\lambda^*\notin\cal S'$. Indeed, since
$\chi_{\lambda}^{2\textrm{-reg}}=\chi_{\lambda^*}^{2\textrm{-reg}}$ (see
Lemma~\ref{ressign}), we
cannot have  $\chi_{\lambda}\in\B$ and $\chi_{\lambda^*}\in\B$
simultaneously
because $\B^{2\textrm{-reg}}$ is free.
Hence, for $\lambda\in\cal S'$, there is no $\lambda'\in\cal S'$
satisfying
$$\cyc{\Ind_{\A_n}^{\sym_n}(\rho_{\lambda}),\Ind_{\A_n}^{\sym_n}
(\rho_{\lambda'})}_{\sym_n}\neq
0.$$
Moreover,
$\Ind_{\A_n}^{\sym_n}(\rho_{\lambda})=\chi_{\lambda}+\chi_{\lambda^*}$
implies
$$\Ind_{\A_n}^{\sym_n}(\rho_{\lambda}^{2\textrm{-reg}})=2\chi_{\lambda}^{2\textrm{-reg}}.$$
(Note that this holds because, for any class function $\alpha$ of $\A_n$,
we have $\Ind_{\A_n}^{\sym_n}(\alpha^{2\textrm{-reg}}) =
(\Ind_{\A_n}^{\sym_n}(\alpha))^{2\textrm{-reg}}$.) Therefore, inducing
Relation~(\ref{eqlibre}) from $\A_n$ to $\sym_n$
we deduce
$$\sum_{\lambda\in\cal S'}2a_{\lambda}\chi_{\lambda}^{2\textrm{-reg}}
+\sum_{\lambda\in\cal S}
(b_{\lambda,-}+b_{\lambda,+})\chi_{\lambda}^{2\textrm{-reg}}=0.$$
But $\B^{2\textrm{-reg}}$ is free, implying $a_{\lambda}=0$ for
$\lambda\in\cal S'$ and $b_{\lambda,+}+b_{\lambda,-}=0$.
Relation~(\ref{eqlibre}) gives
$$\sum_{\lambda\in\cal S}
b_{\lambda,+}(\rho_{\lambda,+}^{2\textrm{-reg}}-\rho_{\lambda,-}^{2\textrm{-reg}})=0.$$
Now, using the fact that $\rho_{\lambda,+}$ and
$\rho_{\lambda,-}$ only differ on the conjugacy classes labelled
by $\overline{\lambda}_{+}$ and $\overline{\lambda}_{-}$
we deduce that $\B_{\A_n}^{2\textrm{-reg}}$ is free (we use here
the same argument as in the proof of~\cite[5.2]{BrGr}).

We now prove that $\B_{\A_n}^{2\textrm{-reg}}$ generates
$\Z\Irr(\A_n)^{2\textrm{-reg}}$ over $\Z$.
Let $\rho$ be a character of $\A_n$ which does not
belong to $\B_{\A_n}$. The definition of $\B_{\A_n}$ implies that there
is $\lambda\vdash n$ with $\lambda\neq\lambda^*$, such that
$\rho=\rho_{\lambda}$. In particular, as explained in~\S\ref{charAn} we
have $$\rho_{\lambda}=\Res_{\A_n}^{\sym_n}(\chi_{\lambda}).$$
Since $\B$ is a $2$-basic set of $\sym_n$, there exist integers
$\{ a_{\chi}, \, \chi \in \B
\}$ satisfying $$\chi_{\lambda}^{\textrm{2-reg}}=\sum_{\chi\in
\B}a_{\chi} \chi^{\textrm{2-reg}}.$$ Restricting this last relation
to $\A_n$, we see 
%The definition of $\B$ implies $\lambda\neq\lambda^*$.
%Hence, $\rho_{\lambda}$ is irreducible.
that
$\rho_{\lambda}^{\textrm{2-reg}}$ is a $\Z$-linear combination of
elements of $B_{\A_n}^{\textrm{2-reg}}$.
This yields the claim.

\section{A $2$-basic set of $\A_n$}\label{part3}

%Throughout this section, for every positive integer $w$, we set
%$$G_w=\Z_2\wr\sym_w.$$

\subsection{A result on wreath products with cyclic
kernel}\label{sectionwr}

Throughout this section, we fix $\ell$ and $w$ positive integers and put
$$G_{\ell,w}=\Z_{\ell}\wr\sym_w,$$
where $\Z_{\ell}$ denotes a cyclic group of order ${\ell}$. We denote by
$\omega$ a generator of $\Z_{\ell}$ and set
$\Irr(\Z_{\ell})=\{\psi_i\,|\,i=1,\ldots,\ell\}$ , with the convention
that $\psi_1$ is the trivial character of $\Z_{\ell}$.
In the following, we denote by $\MP_{\ell,w}$ the
set of $\ell$-tuples of partitions $(\mu_1,\ldots,\mu_{\ell})$ such that
$\sum|\mu_i|=w$.

We recall that the conjugacy classes of $G_{\ell,w}$ are parametrized by
the elements
of $\MP_{\ell,w}$ as follows. The elements of $G_{\ell,w}$ are of the form
$(h,\sigma)$
with $h=(h_1,\ldots,h_w)\in\Z_{\ell}^w$ and $\sigma\in\sym_w$. For any
$k$-cycle
$\kappa=(j,\kappa j,\ldots,\kappa^{k-1} j)$ in $\sigma$, we define
$$g( (h,\sigma);\kappa)=h_jh_{\kappa j}\cdots h_{\kappa^{k-1}
j}\in\Z_{\ell}.$$
Let $\sigma=\prod_{c\in s(\sigma)} c$ be the cycle structure of $\sigma$.
We
then form the corresponding $\ell$-tuples of partitions
$(\mu_1,\ldots,\mu_{\ell})$ by adding
a $k$-cycle to $\mu_i$ whenever $c\in s(\sigma)$ is a
$k$-cycle satisfying $g( (h,\sigma),c)=\omega^{i-1}$.
The resulting $\ell$-tuple $(\mu_1,\ldots,\mu_{\ell})$ lies in
$\MP_{\ell,w}$ and is the so-called
\emph{cycle structure} of $(h,\sigma)$.
Two elements of $G_{\ell,w}$ are conjugate if and only if they have the
same cycle structure.

We define
$$\cal{C}_{\emptyset}=\{(\mu_1,\emptyset,\ldots,\emptyset)\,|\
\mu_1\vdash w\}.$$
\begin{remark}\label{rkinjection}
Using $\sigma\mapsto(1,\sigma)$, we can identify $\sym_w$ to
a subgroup of $G_{\ell,w}$. Note that $\sigma\in\sym_w$ is in the class of
$\sym_w$ labelled by the partition $\mu_1\vdash w$ if and only if
$(1,\mu_1)$ lies in the class of $G_{\ell,w}$ with cycle structure
$(\mu_1,\emptyset,\ldots,\emptyset)\in\cal{C}_{\emptyset}$.
\end{remark}

We recall that the irreducible characters of $G_{\ell,w}$ are also
labelled
by the elements of $\MP_{\ell,w}$ as follows. Let
$\mm=(\mu_1,\ldots,\mu_{\ell})\in\MP_{\ell,w}$.
Consider the character
\begin{equation}
\label{eqphimu}
\phi_{\mm}=\prod_{i=1}^{\ell}\
\underbrace{\psi_i\bt\ldots\bt
\psi_i}_{|\mu_i|\textrm{ times}}.
\end{equation}
If $I_{G_{\ell,w}}(\phi_{\mm})$ denotes the inertial subgroup of
$\phi_{\mm}\in\Irr(\Z_{\ell}^w))$ in $G_{\ell,w}$, then
%$$I_{G_w}(\chi)=\Z_2\wr \sym_{|\mu_1|}\times\Z_2\wr
%\sym_{|\mu_2|}.$$
$$I_{G_{\ell,w}}(\phi_{\mm})=\Z_{\ell}^{w}\rtimes\prod_{i=1}^{\ell}\sym_{|\mu_i|}=\prod_{i=1}^{\ell}\Z_{\ell}\wr
\sym_{|\mu_i|}.$$
Moreover, $\phi_{\mm}$ can be extended to an irreducible character
$\widehat{\phi}_{\mm}=\bt_{i=1}^{\ell}\
\widehat{\psi_i^{|\mu_i|}}$ of
$I_{G_{\ell,w}}(\phi_{\mm})$ by setting
$\widehat{\phi}_{\mm}(h,x)=\phi_{\mm}(h)$
for $h\in\Z_{\ell}^w$
and $x\in\prod\sym_{|\mu_i|}$. The irreducible character $\theta_{\mm}$
corresponding to
$\mm\in\MP_{\ell,w}$ is then given by
\begin{equation}\label{eqtheta}
\theta_{\mm}=\Ind_{I_{G_{\ell,w}}(\phi_{\mm})}^{G_{\ell,w}}(\widehat{\phi_{\mm}}
\otimes
(\chi_{\mu_1}\bt\cdots\bt\chi_{\mu_{\ell}}))=\Ind_{I_{G_{\ell,w}}(\phi_{\mm})}^{G_{\ell,w}}
\left(\prod_{i=1}^{\ell}\
\widehat{\psi_i^{|\mu_i|}} \otimes \chi_{\mu_i}\right),
\end{equation}
where $\chi_{\mu_i}$  denotes the irreducible
characters of $\sym_{|\mu_i|}$ corresponding to the partition $\mu_i$ of
$|\mu_i|$.

\begin{proposition}
\label{valind}
For $\mm=(\mu_1,\ldots,\mu_{\ell})\in\MP_{\ell,w}$, write
$\sym_{\mm}=\sym_{\mu_1}\times\ldots\times\sym_{\mu_{\ell}}$ for the
corresponding Young subgroup of $\sym_w$. Define
$\theta_{\mm}$ as in Formula~(\ref{eqtheta}) and put
$$\Gamma_{\mm}=\Ind_{\sym_{\mm}}^{\sym_w}\left(\chi_{\mu_1}\bt
\cdots\bt \chi_{\mu_{\ell}}\right),$$
where $\chi_{\mu_i}$ denotes the irreducible character of
$\sym_{|\mu_i|}$ corresponding to $\mu_i\vdash|\mu_i|$.
Then, for any $\pi\vdash w$, we have
$$\theta_{\mm}((\pi,\emptyset,\ldots,\emptyset))=\Gamma_{\mm}(\pi).$$
\end{proposition}
\begin{proof}
Let $\pi\vdash w$. We fix $x\in G_{\ell,w}$ in the conjugacy class
of $G_{\ell,w}$ labelled by $(\pi,\emptyset,\ldots,\emptyset)$.
Using~\cite[4.2.10]{James-Kerber}, we deduce
$$|\Cen_{G_{\ell,w}}(x)|=\ell^w\prod_k a_{1k}(x)!k^{a_{1k}(x)},$$
where $a_{1k}(x)$ denotes the number of $k$-parts of the first partition
of the cycle structure of $x$. Denote by $a_k(\pi)$ the number of
$k$-parts of $\pi$. Then we have $a_{1k}(x)=a_{k}(\pi)$
and~\cite[1.2.15]{James-Kerber} implies
\begin{equation}
|\Cen_{G_{\ell,w}}(x)|=\ell^w|\Cen_{\sym_w}(\pi)|.
\label{eqcent}
\end{equation}
Note that
$$\Z_{\ell}^w\rtimes\sym_{\mm}=\prod_{i=1}^{\ell}\Z_{\ell}\wr
\sym_{|\mu_i|}.$$
Then, if we suppose that $\pi=(\pi_1,\ldots,\pi_{l})\in\sym_{\mm}$, it
follows that
$$\Cen_{\Z_{\ell}^w\rtimes\sym_{\mm}}(x)=\prod_{i=1}^{\ell}\Cen_{\Z_{\ell}\wr
\sym_{|\mu_i|}}(\pi_i,\emptyset,\ldots,\emptyset).$$
Furthermore, applying Formula~(\ref{eqcent}) with
$(\pi_i,\emptyset,\ldots,\emptyset)\in G_{\ell,|\mu_i|}$, we deduce
$$
\begin{array}{lcl}
\displaystyle{|\Cen_{\Z_{\ell}^w\rtimes\sym_{\mm}}(x)|}&=&
\displaystyle{\prod_{i=1}^{\ell}\ell^{|\mu_i|}|\Cen_{\sym_{|\mu_i|}}(\pi_i)|}\\
&=&
\displaystyle{\ell^w\prod_{i=1}^{\ell}|\Cen_{\sym_{|\mu_i|}}(\pi_i)|}\\
&=&\ell^w|\Cen_{\sym_{\mm}}(\pi)|.
\end{array}$$
Therefore, the induction formula for characters gives
$$\theta_{\mm}(x)=|\Cen_{G_{\ell,w}}(x)|\sum_{i \in I}
\frac{1}{|\Cen_{\Z_{\ell}^w\rtimes\sym_{\mm}}(x_i)|}\left(\widehat{\phi}_{\mm}
\otimes
(\chi_{\mu_1}\bt\cdots\bt\chi_{\mu_{\ell}})\right)(x_i),$$
where $\widehat{\phi}_{\mm}$ is defined in Equation~(\ref{eqphimu}) and
$\{x_i, \, i \in I\}$ is a system of representatives for the conjugacy
classes of
$\Z_{\ell}^w\rtimes\sym_{\mm}$ such that $x_i$ and $x$ are conjugate in
$G_{\ell,w}$.
However, as explained in Remark~\ref{rkinjection}, the cycle structure of
$x_i$ in
$G_{\ell,w}$ lies in $\cal{C}_{\emptyset}$. Then, for each $i \in I$,
there is $\eta_i\vdash
w$ such that the cycle structure of $x_i$ in $G_{\ell,w}$ is
$(\eta_i,\emptyset,\ldots,\emptyset)$.
%Moreover, $x$ and $x_i$ are conjugate in $G_{\ell,w}$ if and only if
%$\pi$ and $\eta_i$ are conjugate in $\sym_{w}$.
Hence, we have
$$\left(\widehat{\phi}_{\mm}\otimes
(\chi_{\mu_1}\bt\cdots\bt\chi_{\mu_{\ell}})\right)(x_i)=(\chi_{\mu_1}\bt\cdots\bt\chi_{\mu_{\ell}})(\eta_i),$$
because $\widehat{\phi}_{\mm}(1)=\phi_{\mm}(1)=1$.
Moreover,~\cite[4.1]{BrGr} implies that the elements $x_i$ are
conjugate in $\Z_{\ell}^w\rtimes\sym_w$ if and only if the elements
$\eta_i$ are conjugate in $\sym_w$.
We then deduce that the elements $\{ \eta_i, i \in I\}$ form a system of
representatives for the conjugacy classes of $\sym_{\mm}$ such that
$\pi$ and $\eta_i$ are conjugate in $\sym_w$.
It then follows that
$$\begin{array}{lcl}
\theta_{\mm}(x)&=&\displaystyle{\ell^{w}|\Cen_{\sym_w}(\pi)|\sum_{i \in I}
\frac{1}{\ell^w|\Cen_{\sym_{\mm}}(\eta_i)|}
(\chi_{\mu_1}\bt\cdots\bt\chi_{\mu_{\ell}})(\eta_i)}\\
&=&\Gamma_{\mm}(\pi).
\end{array}
$$
\end{proof}

\subsection{A $\Z$-basis of the character ring of
$\sym_{2w}$}\label{part2}
Fix a positive integer $w$. In this section,
%The result of
%this section might be of independent interest.
we will construct a new $\Z$-basis of the character ring of $\sym_{2w}$.

For $\mu\vdash w$, we define
\begin{equation}\label{definduit}
\gamma_{\mu}=\Ind_{\sym_w\times\sym_w}^{\sym_{2w}}(\chi_{\mu}\bt\chi_{\mu}).
\end{equation}
We put
$$B_w=\{\gamma_{\mu}\,|\,\mu\vdash
w\}\cup\{\chi_{\lambda}\,|\,\lambda\notin\cal P_{2w}'\},$$ where
$P_{2w}'$ is the set of partitions defined in Theorem~\ref{main}.

\begin{proposition}
\label{base}
The set $B_w$ is a $\Z$-basis of the ring $\Z\Irr(\sym_{2w})$.
\end{proposition}
\begin{proof}
For $\mu=(\mu_1,\ldots,\mu_r)\vdash w$, we define
$$\widetilde{\mu}=(2\mu_1,\ldots,2\mu_r).$$
Note that the map $\mu\mapsto\widetilde{\mu}$ is a bijection between
${\cal P}_w$ and $\cal P_{2w}'$. It follows that
$|B_w|=|\Irr(\sym_{2w})|$.
Then, to prove that $B_w$ is a $\Z$-basis of $\Z\Irr(\sym_{2w})$, it
is sufficient to prove that $B_w$ is a generating family over $\Z$. Let
$\mu\vdash w$. Then we have
$$\gamma_{\mu}=\sum_{\lambda\vdash
2w}c_{\mu\mu}^{\lambda}\chi_{\lambda},$$
where, for each $\lambda \vdash 2w$, $c_{\mu\mu}^{\lambda}$ is the
Littelwood-Richardson coefficient associated to the partitions $\mu$,
$\mu$ and $\lambda$. If we arrange the partitions of $w$ in
lexicographic order, then~\cite[6.1.2]{GP} implies that the matrix
$$P=(c_{\mu\mu}^{\widetilde{\mu}})_{\mu\vdash w}$$ is a lower triangular
matrix with diagonal entries equal to $1$. In particular, $P$ is
invertible over $\Z$. Then, using $P^{-1}$ we can write the characters
$\chi_{\widetilde{\mu}}$ for $\mu\vdash w$ as a linear combination of
elements of $B_w$ with coefficients in $\Z$. This yields the claim.

\end{proof}

\subsection{A generalized perfect isometry}
First, we will briefly present the notion of generalized perfect
isometry,
introduced in~\cite{KOR} by K\"ulshammer, Olsson and Robinson.
For a union $\cal C$ of conjugacy classes
of a finite group $G$, we say that two irreducible
characters $\alpha$ and $\beta$
are \emph{orthogonal across} $\cal C$ if $$\cyc{\alpha,\beta}_{\cal C}:=
\frac{1}{|G|}\sum_{g\in\cal
C}\alpha(g)\overline{\beta(g)}=0.$$
Then, we define the $\cal C${\emph{-blocks}} of $G$ to be the minimal
non-empty
subsets of $\Irr(G)$ subject to being orthogonal across $\cal C$.
For $b\subseteq\Irr(G)$, we write $(b,\cal C)$ to indicate that $b$ is
a union of $\cal C$-blocks. Then a bijection $\cal I:b\rightarrow b'$
is a {\emph{generalized perfect isometry}} (with respect to ${\cal C}$ and
${\cal C}'$ between two unions of blocks
$(b,\cal C)$ and $(b',\cal{C}')$ of $G$ and $H$, if there are signs
$\{\eta(\alpha)\,|\,\alpha \in b\}$ such that, for all $\alpha,\,\beta\in
b$,
$$\cyc{\alpha,\beta}_{\cal C}=\cyc{\eta(\alpha)\cal{I}(\alpha),
\eta(\beta)\cal{I}(\beta)}_{\cal C'}.$$

Let $w$ be a positive integer. In the following, we put
$$G_w=\Z_2\wr\sym_w.$$
Note that, with the notation of Section~\ref{sectionwr}, we have
$G_w=G_{2,w}$. We can then apply the results of Section~\ref{sectionwr}
to $G_w$. In particular, the irreducible characters of $G_w$ are
labelled by the set $\MP_{2,w}$ and the character corresponding to
$\mm\in\MP_{2,w}$ is denoted by $\theta_{\mm}$.
We also denote by $\cal{C}_{\emptyset}$ the set of elements of $G_w$ with
cycle structure $(\mu_1,\emptyset)$
for some $\mu_1\vdash w$.

%If we denote by
%$\psi_2$ the non-trivial character of $\Z_2$, we remark that
%$\psi_2(-1)=-1$.
In order to describe~\cite[3.6]{BrGr} for $p=2$,
we introduce a bijection on $\MP_{2,w}$ (denoted
by $\check{ \; }$) defined by
$$\check{\mm}=(\mu_1,\mu_2^*)\quad\textrm{for }\mm=(\mu_1,\mu_2).$$

\begin{remark}
This definition comes from \cite{JB}, on which ~\cite[3.6]{BrGr} is
based, and counter-balances the appearance in the Murnaghan-Nakayama
Formula for $G_w$ of some negative signs, coming from the fact that
$\psi_2(-1)=-1$.
\end{remark}

In order to give the main result of this section, we have to recall some
definitions.
For $\lambda\vdash n$, we denote by $\gamma(\lambda)$ the
$2$-core of $\lambda$ (\cite[2.7]{James-Kerber}).
Recall that
two characters $\chi_{\lambda_1}$ and $\chi_{\lambda_2}$ are in
the same $2$-block of $\sym_n$ if and only if the partitions
$\lambda_1$ and $\lambda_2$ have the same $2$-core. Then the
set of $2$-blocks of $\sym_n$ can be labelled by the set of
$2$-cores of $\sym_n$. For $\lambda\vdash n$, the integer
$w=\frac{1}{2}(n-|\gamma(\lambda)|)$ is called the $2$-weight of
$\lambda$.
Since irreducible characters of $\sym_n$ lying in the same $2$-block have
the same $2$-core, it follows that the weight is an invariant of the
block. We can thus define the weight of a block $b$ as the weight of all
characters in $b$.
We now can state:

\begin{theorem}(cf. \cite[3.6]{BrGr}) Let $n$ be a positive integer.
Let $b$ be a $2$-block of $\sym_w$ of $2$-weight $w\neq 0$. For
$\chi_{\lambda}\in b$, write $\alpha_{\lambda}\in\MP_{2,w}$
the $2$-quotient of
$\lambda$. We can associate to $\alpha_{\lambda}$ the character
$\theta_{\alpha_{\lambda}}$ defined in Formula~(\ref{eqphimu}).
Then the map
$$\cal J:b\rightarrow\Irr(G_w),\
\chi_{\lambda}\mapsto\theta_{\check{\alpha_{\lambda}}}$$
is a generalized perfect isometry between $(b,2\textrm{-reg})$ and
$(\Irr(G_w),\cal{C}_{\emptyset})$.
\label{perfectiso}
\end{theorem}
\subsection{Proof of Theorem~\ref{main}}
We keep the notation of the above sections. Recall that,
by~\cite[3.1]{BrGr}, it is possible to define 2-quotients of
partitions in such a way that, for any $n$ and any $\lambda
\vdash n$, if $\alpha_{\lambda}=(\lambda^{(1)},\lambda^{(2)})$,
then $\alpha_{\lambda^*}=(\lambda^{(2)^*},\lambda^{(1)^*})$.
%It
%is well-known that, if $\gamma(\lambda)$ denotes the $2$-core
%of $\lambda$ (\cite[2.7]{James-Kerber}), then $\lambda$ is
%uniquely determined by $\gamma(\lambda)$ and $\alpha_{\lambda}$
%(see~\cite[2.7.30]{James-Kerber}).
Moreover, note that if $\chi_{\lambda}$ lies in a $2$-block of $\sym_n$
with $2$-weight $w$, then $|\lambda^{(1)}|+|\lambda^{(2)}|=w$.

%First, thanks to~\cite[2.1]{BrGr} we can reduce the problem to the
%$2$-blocks of $\sym_n$. Let $b$ be a $2$-block of~$\sym_n$. For every
%$\chi_{\lambda}\in b$ with $2$-quotient $\alpha_{\lambda}=(\mu_1,\mu_2)$,
%the sum $|\mu_1|+|\mu_2|$ is invariant and is called the $2$-weight of $b$.
%In order to prove the theorem, we only have to consider $2$-blocks of $\sym_n$
%with even $2$-weight. Indeed, we have
\begin{lemma}Let $n$ be a positive integer and $\lambda\vdash n$. Then
$\lambda$ is self-conjugate if and only if its 2-quotient has the
form $(\mu,\mu^*)$. In particular, if $w$ is the 2-weight of
$\lambda$, then $w$ has to be even.
\label{2quotient}
\end{lemma}
\begin{proof}Since $\alpha_{\lambda^*}=(\lambda^{(2)^*},\lambda^{(1)^*})$,
it immediately follows that
$\alpha_{\lambda}=(\lambda^{(1)},\lambda^{(1)^*})$ whenever $\lambda\vdash
n$ is a self-conjugate
partition.
\end{proof}
Using this, we can reduce Theorem~\ref{main} to the same question on
$2$-blocks of
$\sym_n$ with even weight. More precisely, we have

\begin{lemma}
The symmetric group $\sym_n$ has a $2$-basic set containing every
$\chi_{\lambda}$ with $\lambda=\lambda^*$ if and only if
every $2$-block $b$ of $\sym_n$ with even weight $w$ has a $2$-basic set
containing
all $\chi_{\lambda}\in b$ with $\lambda=\lambda^*$.
\label{lemred}
\end{lemma}

\begin{proof}
Using~\cite[2.1]{BrGr}, we can reduce the problem to the $2$-blocks of
$\sym_n$. Let $b$ be a $2$-block of $\sym_n$ with odd weight.
Lemma~\ref{2quotient} implies that if $\chi_{\lambda}\in b$, then
$\lambda\neq\lambda^*$. Therefore, it is sufficient to prove that $b$ has
a
$2$-basic set. For this, we use Theorem~\ref{perfectiso}, which implies
that
$\cal J:b\rightarrow\Irr(G_{w})$ is a generalized perfect isometry
between $(b,2\textrm{-reg})$ and $(\Irr(G_{w}),\cal{C}_{\emptyset})$.
Moreover, if we denote  by $B_{\emptyset}$ the set of irreducible
characters of $G_w$ labelled by elements of $\cal{C}_{\emptyset}$, then
Lemma~\cite[4.2]{BrGr} implies $B_{\emptyset}$ is a
$\cal{C}_{\emptyset}$-basic set of $G_w$. The result then follows
from~\cite[2.2]{BrGr}.
\end{proof}

The case of blocks of 2-weight 0 is easy to deal with. Such a block $b$
consists of a unique irreducible character $\chi_{\lambda}$ of $\sym_n$,
such that $\lambda$ is a self-conjugate partition 
(since it is its own 2-core, and one shows easily that any 2-core must be
self-conjugate). Hence $\{ \chi_{\lambda} \}$ is a 2-basic set for $b$, which
obviously satisfies the required property. We next solve the case of
blocks with positive 2-weight.

\begin{proposition}\label{basicset}
Let $b$ be a $2$-block of $\sym_n$ with even weight $2w$ for some
positive integer $w$. Then $b$ has a $2$-basic set containing all
irreducible characters of $b$ labelled by self-conjugate partitions.
%Let $w$ be a non-negative integer and $G_{2w}$ as above. Then
%$G_{2w}$ has a $\cal{C}_{\emptyset}$-basic set containing the
%irreducible characters of $G_w$ labelled by $(\mu,\mu)$.
\end{proposition}

\begin{proof}

Fix a $2$-block $b$ of $\sym_n$ with even $2$-weight $2w$.
Theorem~\ref{perfectiso} implies that
$\cal{J}:b\rightarrow\Irr(G_{2w})$ is a generalized perfect isometry
between $(b,2\textrm{-reg})$ and $(\Irr(G_{2w}),\cal{C}_{\emptyset})$.
We parametrize the irreducible characters of $G_{2w}$ by the elements of
$\MP_{2,2w}$ as described above, and the character corresponding to
$\mm\in\MP_{2,2w}$ will be denoted by $\theta_{\mm}$ as in
Equation~(\ref{eqtheta}).
%Suppose that $b$ has a $2$-basic set containing all irreducible
%characters of $b$ labelled by self-conjugate partition.
Let $\chi_{\lambda}\in b$ with $\lambda=\lambda^*$. Then,
Lemma~\ref{2quotient} and the definition of $\cal J$ imply that
$$\cal J(\chi_{\lambda})=\theta_{(\mu,\mu)}.$$
Therefore, using~\cite[2.2]{BrGr}, we see that proving that $b$ has a
$2$-basic set containing all irreducible characters of $b$ labelled by
self-conjugate partitions is equivalent to showing that the group
$G_{2w}$ has a $\cal{C}_{\emptyset}$-basic set containing all
irreducible characters of $G_{2w}$ labelled by bi-partitions of the form
$(\mu,\mu)$.

Let $B_{\emptyset}$ be the set
of irreducible characters of $G_{2w}$ labelled by elements of
$\cal{C}_{\emptyset}$. More precisely,
the characters of $B_{\emptyset}$ are the characters of $G_{2w}$
with $\Z_2^{2w}$ in their kernel. As we mentioned in Lemma~\ref{lemred},
$B_{\emptyset}$ is a $\cal{C}_{\emptyset}$-basic set of $G_{2w}$.
Note that, for all $\lambda, \, \pi \vdash 2w$, we have
\begin{equation}
\label{eqvalsym}
\theta_{(\lambda,\emptyset)}(\pi,\emptyset)=\chi_{\lambda}(\pi).
\end{equation}
However, $B_{\emptyset}$ doesn't have the required property.
We will now construct from $B_{\emptyset}$ a $\cal{C}_{\emptyset}$-basic
set containing the set of characters $\{\theta_{(\mu,\mu)}\,|\,\mu\vdash
w\}$.
Proposition~\ref{valind} implies that
\begin{equation}
\label{eqvalinduit}
\theta_{(\mu,\mu)}(\pi,\emptyset)=\Gamma_{(\mu,\mu)}(\pi).
\end{equation}
Note that $\Gamma_{(\mu,\mu)}$ is the character $\gamma_{\mu}$
defined in
Formula~(\ref{definduit}).
Furthermore, Formulae~(\ref{eqvalsym}) and~(\ref{eqvalinduit}) imply that,
if we can find a $\Z$-basis of $\sym_{2w}$ containing the character
$\gamma_{\mu}$ for every $\mu\vdash w$ and the irreducible characters
$\chi_{\lambda}$ for $\lambda$ in some parametrizing set $I$, then the set
of irreducible
characters of
$G_w$ labelled by $\{ (\mu,\mu), \, \mu \vdash w \} \cup \{
(\lambda,\emptyset) , \, \lambda\in I \}$ is a
$\cal{C}_{\emptyset}$-basic set of $G_{2w}$. We therefore get the
desired result, with $I=\cal{P}_{2w}\backslash\cal{P}_{2w}'$, by
Proposition~\ref{base}.
\end{proof}
\begin{remark}
Note that the generalized perfect isometry $\cal J$ of
Theorem~\ref{perfectiso} is not one of the isometry described by Osima
between
$b$ and $\Irr(G_w)$. Indeed, Osima's isometry is a generalized perfect
isometry between $(b,2\textrm{-reg})$ and
$(\Irr(G_w),\cal{D}_{\emptyset})$ where $\cal{D}_{\emptyset}$ is the set
of elements of cycle structure
$(\mu_1,\mu_2)\in\MP_{2,w}$ with $\mu_1=\emptyset$.
It seems to be more difficult to prove a result similar to
Proposition~\ref{basicset} for $(\Irr(G_{2w}),\cal{D}_{\emptyset})$.
\end{remark}

We now can prove Theorem~\ref{main}. As explained in the proof of
Lemma~\ref{lemred}, it is sufficient to construct a $2$-basic set with
the required property for all $2$-blocks of $\sym_n$. Let $b$ be a
$2$-block of $\sym_n$ with weight $w$.

If $w$ is odd, then we choose the $2$-basic set $B_b$ of $b$ constructed
in the proof of Lemma~\ref{lemred}. In this case, $\cal{P}_{w}'$ is
empty and the $2$-quotients of the characters in $B_b$ have the form
$(\mu,\emptyset)$ with $\mu\vdash w$.

If $w$ is odd, then Proposition~\ref{basicset} implies that $b$ has a $2$-basic
set $B_b$ satisfying $\chi_{\lambda}\in B_b$ if and only if
$\alpha_{\lambda}=(\mu,\emptyset)$ with
$\mu\in\cal{P}_w\backslash\cal{P}_w'$ or $\alpha_{\lambda}=(\mu,\mu^*)$
for $\mu\vdash w/2$.

Thus, the set $\B_{\Lambda}$ defined in Theorem~\ref{main} is a
$2$-basic set of $\sym_n$.

\bibliographystyle{plain}
\bibliography{references2}
\end{document}